\documentclass{amsart}
\usepackage{amssymb}
\DeclareFontFamily{OMS}{cmsy}{%
\fontdimen16\font=3pt
\fontdimen17\font=3pt}
\def\dj{d\kern-.30em\raise1.25ex\vbox{\hrule width .3em height .03em}}
\def\Dj{D\kern-.70em\raise0.75ex\vbox{\hrule width .3em height .03em}
\kern.03em}
\newenvironment{pf}{\proof[\proofname]}{\endproof}
\newtheorem{lem}{Lemma}
\newtheorem{pro}[lem]{Proposition}
\theoremstyle{definition}
\newtheorem{defn}{Definition}[section]
\def\1{\varnothing}
\def\grten{\mathbin{\widehat{\otimes}}}
\def\adj{\varpi}
\def\inv{{i\!\hspace{0.8pt}n\!\hspace{0.6pt}v}}
\def\bla#1{$(${\it #1\/{}}$)$}
\def\Sum{{\displaystyle\sum}}
\def\e{\epsilon}
\def\k{\kappa}
\def\flO{\Psi}
\def\M{\mathrm{M}}

\def\cal{\mathcal}
\def\Bbb{\mathbb}

\def\id{\mathrm{id}}
\def\BG{\mathrm{B}_G}
\def\g{\mathfrak{g}}
\def\QBG#1{\mathrm{QB}_G^{#1}}
\def\QEG#1{\mathrm{QE}_G^{#1}}
\def\EG{\mathrm{E}_G}
\def\O#1{\cal{O}_{#1}}
\def\pO#1{\cal{H}_{#1}}
\def\OG#1{\cal{O}_{#1,G}}
\def\pOG#1{\cal{H}_{#1,G}}
\def\sOG#1{\cal{Q}_{#1}}
\def\fO{\delta}
\def\mO{\mu}
\def\mM{\mu}
\def\WG{\Omega}
\def\WGi{\daleth}
\def\fWG{\widetilde{\adj}}
\begin{document}
\title[Quantum Principal Bundles]{Quantum Classifying Spaces
$\&$ Universal \\ Quantum Characteristic Classes}
\author{Mi\'co \Dj ur\Dj evi\'c}
\address{Instituto de Matematicas, UNAM, Area de la Investigacion
Cientifica, Circuito Exterior, Ciudad Universitaria, M\'exico DF,
CP 04510, MEXICO}
\thanks{Lectures presented at Quantum Groups $\&$ Quantum
Spaces Minisemester, Stefan Banach International Mathematical
Center, Warsaw, Poland, November 1995.}
\begin{abstract}
A construction of the noncommutative-geometric counterparts of classical
classifying spaces is presented, for general compact matrix quantum
structure groups. A quantum analogue of the classical concept of the
classifying map is introduced and analyzed. Interrelations with the
abstract algebraic theory of quantum characteristic classes are
discussed. Various non-equivalent approaches to defining universal
characteristic classes are outlined. 
\end{abstract}
\maketitle
\section{Introduction}
\renewcommand{\thepage}{}

Let us begin by recalling basic concepts of the
classical classification theory of principal bundles.
From the point of view of this theory, a given
compact Lie group $G$ is represented by the corresponding
{\it classifying space} $\BG$. This space can be equiped with a
(infinite-dimensional) smooth manifold structure.
It classifies principal $G$-bundles in the following sense.
For every compact smooth manifold $M$, isomorphism classes of principal
$G$-bundles $P$ over $M$ are in a natural correspondence with homotopy
classes of maps $\Phi\colon M\rightarrow \BG$. For a given $\Phi$,
the bundle $P$ is reconstructed as the pull-back of a natural contractible
principal $G$-bundle $\EG{}$ over $\BG{}$, via $\Phi$.
This scheme allows us to
associate intrinsically cohomological invariants of $M$, to principal bundles
$P$. Such {\it characteristic classes} are defined as
pull-backs, via $\Phi$, of the cohomology classes of $\BG$. The elements of
the algebra $H(\BG)$ are called {\it universal characteristic classes}.

The algebra $H(\BG)$ can be also described in a completely different way,
without appealing to the construction of the classifying space $\BG$.
Let $\g$ be the Lie algebra of $G$,
and $P(\g)$ the algebra of polynomial functions on $\g$.
Let $\Sigma(\g)\subseteq P(\g)$ be the subalgebra consisting of elements
invariant under the coadjoint action of $G$ on $P(\g)$. Then it can be
shown that
$$H(\BG)=\Sigma(\g),$$
in a natural manner (assuming complex coefficients and the standard
cohomology theory).

In terms of the above identification, the association of
characteristic classes is given via the classical
Weil homomorphism $W\colon\Sigma(\g)\rightarrow H(M)$, obtained
by formally replacing generators of $P(\g)$ by the curvature of an arbitrary
connection on the bundle $P$.

In this letter we are going to discuss how to incorporate
the classical theory of classifying spaces into the conceptual framework of
noncommutative differential geometry \cite{c}. All considerations are
logically based on a general theory of quantum
principal bundles \cite{d1,d2}, in which the bundle and the base are quantum
objects and quantum groups play the role of structure groups.
A detailed exposition of the theory of quantum classifying spaces will
be given in \cite{d4}. 

The first part of the paper
is devoted to the construction of the quantum classifying spaces,
and universal
quantum principal bundles over them. This will be done for general compact
\cite{w1}
matrix quantum groups $G$. Then we shall explain in which sense
such spaces are
`classifying'. The essential point consists in establishing a variant of
the cross-product structuralization \cite{dr-cross} of the
algebra describing an arbitrary
quantum $G$-bundle $P$ over a quantum space $M$. Geometrically, introducing
such a cross-product is equivalent to the specification of
a `classifying map'.

The second part of the paper is devoted to the study of the relations between
cohomology classes associated to a
quantum classifying space in the framework of the general theory
of quantum principal bundles \cite{d2}, and `purely algebraic' universal
characteristic classes analyzed in \cite{d3,d-bc} generalizing the classical
Weil approach. Finally, in the last section some concluding remarks are made.
\renewcommand{\thepage}{\arabic{page}}

\section{The Construction of Quantum Classifying Spaces}

Let $G$ be an arbitrary compact matrix quantum group \cite{w1}. Let $\cal{A}$
be the Hopf *-algebra representing polynomial functions on the
quantum space $G$. We shall denote by
$\phi\colon\cal{A}\rightarrow\cal{A}\otimes\cal{A}$,
$\e\colon\cal{A}\rightarrow\Bbb{C}$ and $\k\colon\cal{A}\rightarrow\cal{A}$
the coproduct, counit and the antipode respectively. 

Let us fix a unitar fundamental representation $u\in\M_n(\cal{A})$ of $G$
such that the conjugate representation $\bar{u}$ is realizable as
a direct summand in $u^m$, for some $m$.
Equivalently, we can say that the matrix elements $u_{ij}$ generate the whole
algebra $\cal{A}$. This is a generalization of the classical unimodularity
condition for $u$. 

For each $d\in\Bbb{N}$, let us denote by $\O{d}$ the unital *-algebra
generated by elements $\psi_{ki}$ where $k\in\{1,...,d\}$ and
$i\in\{1,...,n\}$, and the following relations
\begin{equation}
\sum_k\psi_{ki}^*\psi_{kj}=\delta_{ij}1. 
\end{equation}

By definition, there exists a natural action map $\fO\colon
\O{d}\rightarrow\O{d}\otimes\cal{A}$ defined on the generators by
\begin{equation}
\fO(\psi_{ki})=\sum_j\psi_{kj}\otimes u_{ji}, 
\end{equation}
and extended to the whole $\O{d}$ by requiring that $\delta$ is a
unital *-homomorphism. Let us denote by $\pO{d}$ a unital 
subalgebra of $\O{d}$ generated by the elements $\psi_{ki}$. By definition,
$$ \fO(\pO{d})\subseteq\pO{d}\otimes\cal{A}. $$
We shall denote by $\OG{d}\subseteq\O{d}$ and $\pOG{d}\subseteq\pO{d}$ the
corresponding subalgebras of $G$-invariant elements. 

It is possible to find the $m$-th order elements
$\widehat{\psi}_{\alpha i}\in\pO{d}$
forming the irreducible $\bar{u}$-multiplet, with
\begin{equation*}
\widehat{\psi}_{\alpha i}=\sum_{\omega}
a_{i\omega}\prod_{\begin{subarray}{c}
l\in\alpha\\j\in\omega\end{subarray}}\psi_{lj}
\end{equation*}
where $\alpha$ and $\omega$ are appropriate multiindexes and
$a_{i\omega}\in\Bbb{C}$, so that we have
\begin{gather}
\fO(\widehat{\psi}_{\alpha i})=\sum_j\widehat{\psi}_{\alpha j}
\otimes u_{ji}^*\\
\sum_{\alpha}\widehat{\psi}_{\alpha i}^*\widehat{\psi}_{\alpha j}=
C_{ji}1. 
\end{gather}
Here $C\in \M_n(\Bbb{C})$ is the canonical intertwiner \cite{w1}
between $u$ and the second contragradient representation $u^{\!cc}$. 

The formula
\begin{equation}
S_{\alpha k}=\sum_{ij}\widehat{\psi}_{\alpha i}C_{ji}^{-1}\psi_{kj}
\end{equation}
defines a  matrix $S_{\alpha k}$. It is easy to see that the 
introduced matrix elements are actually from $\pOG{d}$. 

Let us now consider the maps $\tau_{kl}\colon\O{d}\rightarrow\O{d}$, 
defined by 
\begin{equation}
\tau_{kl}(a)=\sum_i\psi_{ki}a\psi_{li}^*. 
\end{equation}
We have then
\begin{equation}
\sum_r\tau_{kr}(a)\tau_{rl}(b)=\tau_{kl}(ab)\qquad\tau_{kl}(a)^*
=\tau_{lk}(a^*), 
\end{equation}
in other words the associated
matrix map $\tau\colon\O{d}\rightarrow\M_d(\O{d})$ is a *-homomorphism. 
From the definition of $\tau$ it follows that the algebra 
$\OG{d}$ is $\tau$-invariant. In what follows, the domain of $\tau$ will be
restricted to $\OG{d}$. 

The algebra $\O{d}$ can be naturally decomposed as
$$
\O{d}=\sideset{}{^\oplus}\sum_{r\in\cal{T}}\O{d}^r
$$
where $\cal{T}$ is a complete set of mutually inequivalent irreducible
representations of $G$, and
$\O{d}^r\subseteq\O{d}$ are the multiple irreducible $\OG{d}$-submodules
corresponding to the decomposition of $\fO$. Obviously we have
$\O{d}^\1=\OG{d}$, where $\1$ is the trivial representation of $G$. 
We can similarly decompose
the $\delta$-invariant subalgebra $\pO{d}$.

Let us denote by $\sOG{d}\subseteq\OG{d}$ the minimal
unital $\tau$-invariant subalgebra containing all
the elements $S_{\alpha k}^*$.

\begin{lem} \bla{i} The following commutation relations hold
\begin{align}
\psi_{ki}a&=\sum_l\tau_{kl}(a)\psi_{li}\label{psi-a}\\
\psi_{ki}^*a&=\sum_{\alpha\beta}S_{\alpha k}^*\tau^m_{\alpha\beta}(a)
\widehat{\psi}_{\beta i}\label{psi*-a}. 
\end{align}

\bla{ii}
Every element $\psi$ of a multiple irreducible module 
$\O{d}^r$ can be written in the form 
\begin{equation}\label{dec-r}
\psi=\sum_{f\varphi} f\varphi,
\end{equation}
where $\varphi\in\pO{d}^r$ and $f\in\sOG{d}$. 
\end{lem}
\begin{pf}
The statement \bla{i} follows from definitions of $S_{\alpha k}$, $\tau$
and $\widehat{\psi}_{\alpha i}$. Decomposition~\eqref{dec-r} follows from
commutation relations~\eqref{psi-a} and~\eqref{psi*-a}. 
\end{pf}
\section{Reconstruction $\&$ Classifying Maps}

Let us consider a unital *-algebra $\cal{V}$, interpreted as consisting of
`smooth functions' on a quantum space $M$. Let $P=(\cal{B},i,F)$ be a
quantum principal $G$-bundle \cite{d2} over $M$. Here $\cal{B}$ is a *-algebra
representing $P$ as a quantum space, while $F\colon\cal{B}\rightarrow\cal{B}
\otimes\cal{A}$ and $i\colon\cal{V}\rightarrow\cal{B}$ are *-homomorphisms
corresponding to the right action of $G$ on $P$, and to the projection of $P$
onto $M$ respectively.

Explicitly, the following properties hold
\begin{gather*}
(\id\otimes\phi)F=(F\otimes\id)F\qquad (\id\otimes\e)F=\id\\
i(\cal{V})=\Bigl\{\,b\in\cal{B}\mid F(b)=b\otimes1\,\Bigr\},
\end{gather*}
and for each $a\in\cal{A}$ there exist elements
$q_\alpha,b_\alpha\in\cal{B}$ satisfying
$$\sum_\alpha b_\alpha F(q_\alpha)=1\otimes a.$$
The last condition corresponds to the classical requirement that the
structure group is acting freely on the bundle. 

It is possible to prove \cite{d4} that under relatively weak `regularity'
conditions for the bundle $P$, there exist a natural number $d$ and a
$\cal{B}$-matrix
$B=\{b_{ki}\}$,
where $k\in\{1,\dots,d\}$ and $i\in\{1,\dots,n\}$ such that the following
relations hold:
\begin{gather}
F(b_{ki})=\sum_j b_{kj}\otimes u_{ji},\\
\sum_k b_{ki}^* b_{kj}=\delta_{ij}1. 
\end{gather}

\begin{defn} We say that a quantum principal bundle
$P$ admitting a matrix $B$ of the above described type has a
complexity level $d$. 
\end{defn}

Let us observe that the number $d$ figuring in the above definition is
not fixed uniquely-we can always increase this number as necessary.
Having a matrix $B$ we can introduce a homomorphism $\rho\colon
\cal{V}\rightarrow \M_d(\cal{V})$, by the formula
\begin{equation}
\rho_{kl}(f)=\sum_i b_{ki}f b_{li}^*. 
\end{equation}
It follows immediately that the following commutation relations hold
\begin{equation}
b_{ki}f=\sum_l\rho_{kl}(f)b_{li},
\end{equation}
between the elements of $B$ and $\cal{V}$. 

Let $\QBG{d}$ be the quantum space corresponding tp $\OG{d}$, and
let $\iota\colon\OG{d}\rightarrow\O{d}$ be the canonical inclusion map.

\begin{lem} The triplet $\QEG{d}=(\O{d},\iota,\delta)$ is a quantum principal
bundle over $\QBG{d}$, having the complexity level $d$. \qed
\end{lem}

Furthermore, there exists the
unique unital *-homomorphism $\gamma\colon\O{d}\rightarrow\cal{B}$, specified
by
\begin{equation}
\gamma(\psi_{ki})=b_{ki},
\end{equation}
and we have $\gamma(\OG{d})\subseteq\cal{V}$, because $\gamma$ intertwines
$\fO$ and $F$. 

\begin{lem} The following natural decomposition holds
\begin{equation}
\cal{B}\leftrightarrow \cal{V}\otimes_*\O{d}\leftrightarrow
\O{d}\otimes_*\cal{V},
\end{equation}
where $\otimes_*$ denotes the tensor product over $\OG{d}$. \qed
\end{lem}

We are going to reverse the above construction. Let
$\rho\colon\cal{V}\rightarrow\M_d(\cal{V})$ be a (non-trivial) 
*-homomorphism and let $\gamma\colon\OG{d}\rightarrow\cal{V}$ 
be a unital *-homomorphism such that
\begin{gather}
\sum_\beta\rho_{\alpha\beta}(f)\gamma(o_{\beta j})
=\gamma(o_{\alpha j})f\label{ittr}\\
\gamma\tau_{kl}=\rho_{kl}\gamma,
\end{gather}
for every $f\in\cal{V}$ and every multi-indexed system
$$o_{\alpha j}=\sum_\omega c_{\alpha\omega}^j\prod_{\begin{subarray}{c}
k\in\alpha\\i\in\omega\end{subarray}}\psi_{ki}$$
from the algebra $\pOG{d}$.
In formula \eqref{ittr} we have considered the appropriate
iterations of $\rho$. 

Motivated by the results of the above analysis, we shall define a crossproduct
between $\O{d}$ and $\cal{V}$ and reconstruct in such a way the bundle $P$. 
This construction is a generalization of the standard crossproduct
construction considered in \cite{dr-cross}. 

The algebra $\cal{V}$ is a bimodule over $\OG{d}$, in a natural manner. 
The bimodule structure is induced by the map $\gamma$. We shall denote
by $\otimes_\gamma$ the corresponding bimodule tensor products between
algebras $\cal{V}$ and $\O{d}$. 

The essential part of the construction is to prove the existence of the
canonical flip-over map $\flO\colon\O{d}\otimes_\gamma\cal{V}\rightarrow
\cal{V}\otimes_\gamma\O{d}$. The consistency of the whole construction
will be ensured by the compatibility relations between $\gamma$ and $\rho$. 

At first, let us forget for a moment about the relations in $\O{d}$, and 
consider it as a free algebra. Let $\flO\colon\O{d}\otimes\cal{V}
\rightarrow\cal{V}\otimes\O{d}$ be the linear map specified by 
\begin{gather*}
\flO(\psi_{ki}\otimes f)=\sum_l\rho_{kl}(f)\otimes\psi_{li}
\qquad\flO(1\otimes f)=f\otimes 1\\
\flO(\psi_{ki}^*\otimes f)=\sum_{\alpha\beta}R^*_{\alpha k}
\rho^m_{\alpha\beta}(f)\otimes\widehat{\psi}_{\beta i}\\
\flO(\mO\otimes\id)=(\id\otimes \mO)(\flO\otimes\id)(\id\otimes\flO),
\end{gather*}
where $\mO\colon\O{d}\otimes\O{d}\rightarrow\O{d}$ is the 
corresponding product map and we have put
$R_{\alpha k}=\gamma(S_{\alpha k})$.
Let us denote by the same symbol $\mM$ be the product in $\cal{V}$.
\begin{lem} We have
\begin{equation}
\flO(\id\otimes\mM)=(\mM\otimes\id)(\id\otimes\flO)(\flO\otimes\id). 
\end{equation}
\end{lem}
\begin{pf} It is sufficient to check the equality on elements of the form
$\psi_{ki}\otimes f$ and $\psi_{ki}^*\otimes f$. Direct transformations give 
\begin{equation*}
\begin{split}
\flO(\psi_{ki}\otimes fg)&=\sum_l\rho(fg)_{kl}
\otimes\psi_{li}=\sum_{lr}\rho(f)_{kr}\rho(g)_{rl}\otimes\psi_{li}\\
&=\sum_{lr}(\mM\otimes\id)\Bigl\{
\rho(f)_{kr}\otimes\rho(g)_{rl}\otimes\psi_{li}\Bigr\}\\
&=(\mM\otimes\id)(\id\otimes\flO)(\flO\otimes\id)(\psi_{ki}\otimes f
\otimes g)
\end{split}
\end{equation*}
and similarly
\begin{equation*}
\begin{split}
\flO(\psi_{ki}^*\otimes fg)&=\sum_{\alpha\beta}R^*_{\alpha k}
\rho^m_{\alpha\beta}(fg)\otimes\widehat{\psi}_{\beta i}=
\sum_{\alpha\beta\gamma}R^*_{\alpha k}\rho^m(f)_{\alpha\gamma}
\rho^m(g)_{\gamma\beta}\otimes\widehat{\psi}_{\beta i}\\&=
\sum_{\alpha\beta\gamma}(\mM\otimes\id)\Bigl\{
R^*_{\alpha k}\rho^m_{\alpha\gamma}(f)\otimes
\rho^m_{\gamma\beta}(g)\otimes\widehat{\psi}_{\beta i}\Bigr\}\\&=
(\mM\otimes\id)(\id\otimes\flO)(\flO\otimes\id)(\psi_{ki}^*\otimes f
\otimes g).\qed
\end{split}
\end{equation*}
\renewcommand{\qed}{}
\end{pf}

The map $\flO$ preserves the transformation properties of the elements
from $\O{d}$, in a natural manner. 
Let us now assume that the map $\flO$ is `partially factorized' in
the form $\flO
\colon\O{d}\otimes\cal{V}\rightarrow\cal{V}\otimes_\gamma\O{d}$.
This allows us to consider the initial relations. 

\begin{lem} At the level of the partial factorization, it is possible 
to include consistently the defining relations for $\O{d}$ in the game. 
\end{lem}
\begin{pf} We have to check that the defining relations `pass' through the
twist. Applying the definition of $\flO$, commutation 
relations between $\psi_{ki}$ and $\cal{V}$,
 and the fact that the resulting product is
$\gamma$-relativized, we obtain
\begin{equation*}
\begin{split} \flO(\delta_{ij}1\otimes f)=
\sum_k\flO(\psi_{ki}^*\psi_{kj}\otimes f)&=\sum_{kl}(\id\otimes\mO)
(\flO\otimes\id)(\psi_{ki}^*\otimes\rho(f)_{kl}\otimes\psi_{lj})\\
&=\sum_{kl\alpha\beta}R_{\alpha k}^*\rho^m_{\alpha\beta}
\bigl[\rho(f)_{kl}\bigr]\otimes\widehat{\psi}_{\beta i}\psi_{lj}\\
&=\sum_{l\beta}fR_{\beta l}^*\otimes\widehat{\psi}_{\beta i}\psi_{lj}=
\sum_{l\beta}f\otimes S_{\beta l}^*\widehat{\psi}_{\beta i}\psi_{lj}\\
&=f\otimes\delta_{ij}1, 
\end{split}
\end{equation*}
which completes the proof. 
\end{pf}

Let us now assume that the following identity holds
\begin{equation}
\gamma(\psi)f=\sum_kf_k\gamma(\psi_k), 
\end{equation}
where $f\in\cal{V}$ and $\psi\in\sOG{d}$, while $\Sum_kf_k\otimes\psi_k=
\Psi(\psi\otimes f)$. 
\begin{defn}
Every pair $(\gamma,\rho)$ satisfying all the above mentioned conditions
is called a $d$-classifying map for $M$. 
\end{defn}

Let us consider the question of the projectability of twists
to the tensor product
over $\gamma$. Consistent factorization to the tensor products over $\gamma$
will be possible only if the twist $\flO$ is $\OG{d}$-linear, on both sides. 
\begin{lem}
The map $\flO\colon\O{d}\otimes\cal{V}\rightarrow\cal{V}\otimes_\gamma
\O{d}$ is a homomorphism of $\OG{d}$-bimodules. 
\end{lem}
\begin{pf} Let us consider an arbitrary element $\psi\in\OG{d}$. We have
then
\begin{multline*} \flO(\psi x\otimes f)=\sum_k(\id\otimes\mO)\bigl[
\flO(\psi\otimes f_k)\otimes x_k\bigr]=\sum_{kl}f_{kl}\otimes\psi_l x_k=
\sum_{kl}f_{kl}\gamma(\psi_l)\otimes x_k\\
=\sum_k \gamma(\psi)f_k\otimes x_k=\gamma(\psi)\flO(x\otimes f), 
\end{multline*}
where $\flO(x\otimes f)=\Sum_kf_k\otimes x_k$ and
$\flO(\psi\otimes f_k)=\Sum_l f_{kl}\otimes\psi_l$. Now, to prove the 
right $\OG{d}$-linearity, it is sufficient to consider homogeneous 
blocks  $x_{ki}$ from the algebra $\pO{d}$, where $k,i$ are appropriate
multiindexes. A direct computation gives
\begin{equation*}
\begin{split}
\flO\bigl(x_{ki}\otimes f\gamma(\psi)\bigr)&=\sum_l
(\mM\otimes\id)\bigl[ \rho^r_{kl}(f)\otimes
\flO\bigl(x_{li}\otimes\gamma(\psi)\bigr)\bigr]\\&=
\sum_{ln}\rho^r_{kl}(f)(\rho^r_{
ln}\gamma)(\psi)\otimes x_{ni}\\&=\sum_{ln}\rho^r_{kl}(f)\gamma[\tau^r_{ln}(
\psi)]\otimes x_{ni}
=\sum_{ln}\rho^r_{kl}(f)\otimes[\tau^r_{ln}(\psi)]x_{ni}\\&=
\sum_l\rho^r_{kl}(f)\otimes x_{li}\psi=\flO(x_{ki}\otimes f)\psi, 
\end{split}
\end{equation*}
where $r$ is the degree of elements $x_{ki}$. Hence, $\flO$ is 
$\OG{d}$-linear on both sides. 
\end{pf}

Let $\cal{Q}\subseteq\OG{d}$ be the space of elements $\psi$ satisfying
\begin{equation}\label{def-Q}
\flO\bigl(x\psi\otimes f\bigr)=\flO\bigl(x\otimes\gamma(\psi)f\bigr)
\end{equation}
for each $x\in\O{d}$ and $f\in\cal{V}$. By definition, $\cal{Q}$ is a
subalgebra of $\OG{d}$ and $1\in\cal{Q}$. 

\begin{lem} \bla{i} The algebra $\cal{Q}$ contains the elements 
$S_{\alpha k}^*$ and we have $\pOG{d}\subseteq\cal{Q}$.

\smallskip
\bla{ii} The following equality holds
\begin{equation}
\sum_r\flO\bigl(p_{kr}\otimes\rho_{rl}(f)\bigr)=\rho_{kl}(f)\otimes 1, 
\end{equation}
for each $f\in\cal{V}$. Here $p_{kr}=\tau_{kr}(1)$ are matrix elements of
the associated projection in $\M_d(\OG{d})$. 
\end{lem}
\begin{pf} We shall prove here the second statement. A direct computation
gives
\begin{equation*}
\begin{split}
\sum_r\flO\bigl(p_{kr}\otimes\rho_{rl}(f)\bigr)&=\sum_{r\alpha\beta i}
(\id\otimes\mO)(\flO\otimes\id)\Bigl(\psi_{ki}\otimes R^*_{\alpha r}
\rho^m_{\alpha\beta}[\rho_{rl}(f)]\otimes\widehat{\psi}_{\beta i}\Bigr)\\
&=\sum_{rn\alpha\beta i}\rho_{kn}\Bigl\{R_{\alpha r}^*\rho^m_{\alpha\beta}
\bigl[\rho_{rl}(f)\bigr]\Bigr\}\otimes\psi_{ni}   \widehat{\psi}_{\beta i}
\\&=\sum_{n\beta i}\rho_{kn}[fR^*_{\beta l}]\otimes\psi_{ni}\widehat{\psi}_{
\beta i}\\&=\sum_{nr\beta i}\rho_{kr}(f)\otimes\tau_{rn}[S_{\beta l}^*]
\psi_{ni}\widehat{\psi}_{\beta i}=\sum_{r\beta i}\rho_{kr}(f)
\otimes\psi_{ri}S_{\beta l}^*\widehat{\psi}_{\beta i}\\&=
\sum_{ri}\rho_{kr}(f)\otimes\psi_{ri}\psi^*_{li}=\sum_r
\rho_{kr}(f)\otimes p_{rl}=\rho_{kl}(f)\otimes 1. 
\end{split}
\end{equation*}
We have used the definition of the tensor product over $\gamma$, and the
basic commutation relations in $\cal{V}$ and $\O{d}$. 
\end{pf}

\begin{lem} The algebra $\cal{Q}$ is closed under the action of operators
$\tau_{kl}$. 
\end{lem}
\begin{pf} Let us consider an arbitrary element $\varphi\in\cal{Q}$. We have
then
\begin{equation*}
\begin{split}
\flO\bigl(x\tau_{kl}(\varphi)\otimes f\bigr)&=
\sum_{\alpha\beta i}(\id\otimes\mO)
(\flO\otimes\id)\Bigl(x\psi_{ki}\varphi\otimes
R^*_{\alpha l}\rho^m_{\alpha\beta}
(f)\otimes\widehat{\psi}_{\beta i}\Bigr)\\
&=\sum_{\alpha\beta i}(\id\otimes\mO)
(\flO\otimes\id)\Bigl(x\psi_{ki}\otimes\gamma(\varphi)
R^*_{\alpha l}\rho^m_{\alpha\beta}
(f)\otimes\widehat{\psi}_{\beta i}\Bigr)\\
&=\sum_{\alpha\beta ni}(\id\otimes\mO)
(\flO\otimes\id)\Bigl(
x\psi_{ki}\otimes R^*_{\alpha n}\rho^m_{\alpha\beta}
\bigl[(\rho_{nl}\gamma)(\varphi)f\bigr]
\otimes\widehat{\psi}_{\beta i}\Bigr)\\
&=\sum_{ni}(\id\otimes\mO)(\flO\otimes\id)(\id\otimes\flO)
\bigl(x\psi_{ki}
\otimes\psi_{ni}^*\otimes(\rho_{nl}\gamma)(\varphi)f\bigr)\\
&=\sum_r\flO\bigl(xp_{kr}\otimes\rho_{rl}\gamma(\varphi)f\bigr)=
\flO\bigl(x\otimes\rho_{kl}\gamma(\varphi)f\bigr). 
\end{split}
\end{equation*}
The last equality in the above sequence of 
transformations is justified by the previous lemma. 
\end{pf}

From the sequence of the previous lemmas, it follows that
\begin{equation}
\cal{Q}=\OG{d}. 
\end{equation}
In other words, the domain of the map $\flO$ can be factorized 
through $\gamma$. In such a way we obtain the fully factorized 
$\OG{d}$-bimodule map
\begin{equation}
\flO\colon\O{d}\otimes_\gamma\cal{V}\rightarrow\cal{V}\otimes_\gamma
\O{d}.  
\end{equation}
The properties of this map are summarized in the following 

\begin{pro} \bla{i} The following identities hold
\begin{gather}
\flO(\id\otimes\mM)=(\mM\otimes\id)(\id\otimes\flO)(\flO\otimes\id)\\
\flO(\mO\otimes\id)=(\id\otimes\mO)(\flO\otimes\id)(\id\otimes\flO)\\
\flO(\psi_{ki}\otimes f)=\sum_l\rho_{kl}(f)\otimes\psi_{li}\\
\flO(\psi_{ki}^*\otimes f)=\sum_{\alpha\beta}R^*_{\alpha k}
\rho^m_{\alpha\beta}(f)\otimes\widehat{\psi}_{\beta i}. 
\end{gather}

\bla{ii} The map $\flO$ is bijective, and 
\begin{equation}\label{*-flO}
\flO^{-1}=*\flO *
\end{equation}
where $*\colon\cal{V}\otimes_\gamma\O{d}\leftrightarrow
\O{d}\otimes_\gamma\cal{V}$ is the canonical conjugation map.
\end{pro}

\begin{pf} We have to prove property \bla{ii}. Let us assume that 
the element $\psi\in\O{d}$ possesses the property that equality $\flO*\flO
=*$ holds on all elements of the form $\psi\otimes f$. Evidently, 
the space $\cal{L}$ of elements $\psi$ contains the unit $1\in\O{d}$ and it 
is a left $\OG{d}$-submodule of $\O{d}$. 
It turns out that $\cal{L}$ is also a right $\pO{d}$-submodule of $\O{d}$. 
Indeed, 
\begin{equation*}
\begin{split}
(\flO*\flO)(\psi\psi_{ki}\otimes f)&
=\sum_l(\flO*)(\id\otimes\mO)(\flO\otimes\id)
\bigl(\psi\otimes\rho_{kl}(f)\otimes\psi_{li}\bigr)\\&=
\sum_l\flO(\mO\otimes\id)\bigl[\psi^*_{li}\otimes(*\flO)\bigl(\psi
\otimes\rho_{kl}(f)\bigr)\bigr]\\
&=\sum_l(\id\otimes\mO)(\flO\otimes\id)\bigl[\psi^*_{li}\otimes
(\flO*\flO)\bigl(\psi
\otimes\rho_{kl}(f)\bigr)\bigr]\\&=
\sum_l(\id\otimes\mO)(\flO\otimes\id)(\psi_{li}^*
\otimes\rho_{lk}(f)^*\otimes\psi^*)\\
&=\sum_{l\alpha\beta}R^*_{\alpha l}\rho^m_{\alpha\beta}\rho_{lk}(f^*)
\otimes\widehat{\psi}_{\beta i}\psi^*=\sum_\beta 
f^*R^*_{\beta k}\otimes\widehat{\psi}_{\beta i}\psi^*\\
&=\sum_\beta f^*\otimes S^*_{\beta k}\widehat{\psi}_{\beta i}\psi^*
=f^*\otimes(\psi\psi_{ki})^*. 
\end{split}
\end{equation*}
Hence, $\cal{L}=\O{d}$, which completes the proof. 
\end{pf}

As a direct consequence of the above proposition, we see that formulas
\begin{align}
(g\otimes\psi)(f\otimes\varphi)&=g\flO(\psi\otimes f)\varphi\\
(f\otimes\psi)^*&=\flO(\psi^*\otimes f^*)
\end{align}
consistently define a *-algebra structure on the space
$$\cal{B}=\cal{V}\otimes_\gamma\O{d}. $$
Furthermore, a natural action $\fO$ of $G$ on $\O{d}$ 
is naturally projectable to a comodule *-algebra map 
$F\colon\cal{B}\rightarrow\cal{B}\otimes\cal{A}$. Let 
$i\colon\cal{V}\rightarrow\cal{B}$ be the canonical inclusion map, 
given by $i(f)=f\otimes 1$. By construction, $i$ is a unital *-homomorphism. 

\begin{pro} The triplet $P=(\cal{B},i,F)$ is a quantum principal $G$-bundle
over the quantum space $M$. 
\end{pro}

\begin{pf} Let us consider the vertical integration map
$h_M\colon\cal{B}\rightarrow\cal{B}$ given by 
$h_M=(\id\otimes h)F$, where $h\colon\cal{A}\rightarrow\Bbb{C}$ is the Haar
measure \cite{w1} of $G$. We have
$$ h_Mi(f)=f\qquad h_M(\cal{B})=\cal{V},$$
as follows from the definitions of $\cal{B}$ and $F$. 
Hence $i$ is injective and its image coincides with the $F$-fixed point 
subalgebra of $\cal{B}$. The freeness of $F$ is a direct consequence of the
freeness of $\fO$. 
\end{pf}

Therefore, we have the following natural correspondence
\begin{equation*}
\left\{\!
\begin{gathered}
\mbox{Quantum principal
$G$-bundles $P$}\\
\mbox{having complexity level $d$}\\
\mbox{equipped with a $d\times d$ matrix $B$}
\end{gathered}\right\}\leftrightarrow\Bigl\{\mbox{$d$-classifying
maps $(\rho,\gamma)$}\Bigr\}.
\end{equation*}
The above correspondence simplifies if we factorize down to the appropriate
homotopy classes. Namely, it turns out \cite{d4} that the homotopy
class of $B$, as well as the homotopy class of the
classifying map $(\rho,\gamma)$ are stable, if we stay in the
framework of the same homotopy class of $P$
and allow $d$ to take arbitrarily high values (if necessary).
Factorizing through the homotopy equivalence we obtain the correspondence
\begin{equation*}
\left\{\!
\begin{gathered}
\mbox{Homotopy classes of quantum principal}\\
\mbox{$G$-bundles $P$ over $M$}\\
\mbox{having finite complexity levels}
\end{gathered}\right\}\leftrightarrow
\left\{
\begin{gathered}
\mbox{Homotopy classes}\\
\mbox{of classifying maps $(\rho,\gamma)$}
\end{gathered}
\right\}.
\end{equation*}

\section{Quantum Characteristic Classes}

In this section we shall consider a quantum generalization
\cite{d3,d-bc} of the classical Weil
homomorphism, and discuss its relations with the quantum classifying spaces. 

Let $\Gamma$ be a bicovariant \cite{w2} first-order *-calculus over $G$. 
Let us consider the universal graded-differential
algebra $\WG$ built over the vector space $\Gamma_{\inv}$ of left-invariant
elements of $\Gamma$. We assume that
$d(1)=0$, and that the elements of $\Gamma_{\inv}$ have degree one in $\WG$. 
The natural *-involution on $\Gamma_{\inv}$ naturally extends to $\WG$,
so that the
differential $d\colon\WG\rightarrow\WG$ is hermitian. By construction,
we have $H(\WG)=\Bbb{C}$.

To fix ideas, we shall assume that the higher-order calculus on $G$ is
given by the universal \cite{d1} differential envelope $\Gamma^\wedge$ of
$\Gamma$. Another natural choice is to assume that
the higher-order calculus is based on the
braided exterior \cite{w2} algebra $\Gamma^\vee$, associated to
$\Gamma$. 

Let $\fWG\colon\WG\rightarrow\WG\grten\Gamma^\wedge$ be the unital
graded-differential homomorphism uniquely characterized by
the property
$$
\fWG(\vartheta)=\adj(\vartheta)+1\otimes\vartheta
$$
for each $\vartheta\in\Gamma_{\inv}$. Here
$\adj\colon\Gamma_{\inv}\rightarrow\Gamma_{\inv}\otimes\cal{A}$ is
the corresponding adjoint action (the restriction to
$\Gamma_{\inv}$ of the right action $\varrho_\Gamma\colon
\Gamma\rightarrow\Gamma\otimes\cal{A}$).
The map $\fWG$ is hermitian, and we have
$$
(\fWG\otimes\id)\fWG=(\id\otimes\widehat{\phi})\fWG,
$$
where $\widehat{\phi}\colon\Gamma^\wedge\rightarrow\Gamma^\wedge\grten
\Gamma^\wedge$ is the graded-differential extension of the coproduct map
(acting on $\Gamma$ as a direct sum of left/right action maps
$\ell_\Gamma\colon\Gamma\rightarrow
\cal{A}\otimes\Gamma$ and $\varrho_\Gamma$).

Let us now consider the subalgebra
$\WGi\subseteq\WG$ consisting of all the elements invariant under
the action $\fWG$. This algebra generally possesses non-trivial
cohomology classes.

\begin{defn}\label{def:uni} The elements of the graded *-algebra
$H(\WGi)$ are called
{\it universal characteristic classes} associated to $G$ and $\Gamma$. 
\end{defn}

Let us consider an
arbitrary $P=(\cal{B},i,F)$ quantum principal $G$-bundle over $M$,
and let $\Omega(P)$ be an arbitrary differential calculus over $P$,
in the sense of \cite{d2}. This means that $\Omega(P)$ is a
graded-differential *-algebra extending $\cal{B}$, generated by $\cal{B}$,
and such that there exists a graded-differential extension of
$\widehat{F}\colon\Omega(P)\rightarrow\Omega(P)\grten\Gamma^\wedge$ of the
right action $F$. Let $\Omega(M)\subseteq\Omega(P)$ be a graded-differential
*-subalgebra consisting of $\widehat{F}$-invariant elements. The elements
of this subalgebra play the role of differential forms on the base space $M$. 

Finally, let us consider an arbitrary connection $\omega$ on $P$.
By definition \cite{d2}, this means that $\omega\colon
\Gamma_{\inv}\rightarrow\Omega(P)$ is a first-order hermitian map satisfying
$$ \widehat{F}\omega(\vartheta)=\adj(\vartheta)+1\otimes\vartheta$$
for each $\vartheta\in\Gamma_{\inv}$. By universality of the algebra $\Omega$,
this map admits the unique extension
$\widehat{\omega}\colon\WG\rightarrow\Omega(P)$, which is a
graded-differential algebra homomorphism. Furthermore, $\widehat{\omega}$ is
hermitian, and it intertwines maps $\fWG$ and
$\widehat{F}$. This intertwining property implies that
$\widehat{\omega}(\WGi)\subseteq\Omega(M)$.

\begin{lem} The induced cohomology map $W\colon H(\WGi)\rightarrow
H(M)$ is independent of the choice of $\omega$. \qed
\end{lem}

The above lemma justifies the definition~\ref{def:uni}. The main idea is
to define characteristic classes as generic cohomology classes expressed
algebraically via the connection $\omega$ and its differential $d\omega$. 
The map $W$ is
a general quantum counterpart of the Weil homomorphism.

Let us now apply the construction to the universal bundle,
$\QEG{d}$ over the classifying space $\QBG{d}$, assuming that the calculus
on $\QEG{d}$ is based on the appropriate universal envelope. 
Let $W_G\colon H(\WGi)\rightarrow H(\QBG{d})$
be the corresponding Weil homomorphism. 

\begin{pro} The map $W_G$ is bijective.
Therefore, we can naturally identify cohomology algebras $H(\WGi)$ and
$H(\QBG{d})$.\qed
\end{pro}

On the other hand, for an arbitrary complexity-$d$ bundle $P$ with a
differential structure $\Omega(P)$, the natural map
$\gamma\colon\O{d}\rightarrow\cal{B}$ uniquely extends to a differential
algebra homomorphism $\gamma\colon\Omega(\QEG{d})\rightarrow
\Omega(P)$, which satisfies $\gamma\bigl[\Omega(\QBG{d})\bigr]
\subseteq\Omega(M)$. This is a consequence of the fact that $\gamma$
intertwines $\widehat{\delta}$ and $\widehat{F}$.

Let $\gamma_*\colon H(\QBG{d})\rightarrow H(M)$ be
the corresponding cohomology map.

\begin{lem} Under the above assumptions we have
\begin{equation}
W=\gamma_* W_G. 
\end{equation}
In other words, the quantum Weil homomorphism is induced by an arbitrary
classifying map. \qed
\end{lem}

For the end of this section, let us mention that there exists
another \cite{d3,d-bc},
nonequivalent, approach to defining universal characteristic classes,
which is from the conceptual point of view closer to the original Weil
construction. This approach is applicable only to bundles equipped with
differential structures admitting very special connections (called
regular and multiplicative \cite{d2}). In this case it is possible
to construct, with the help of the curvature of an arbitrary
regular and multiplicative connection, a natural map
$w\colon\Sigma\rightarrow H\!Z(M)$, where $H\!Z(M)$ is the cohomology
of the graded centre of $\Omega(M)$ and $\Sigma$ is the $G$-invariant part
of the braided-symmetric algebra built over $\Gamma_{\inv}$, relative to
the natural \cite{w2} braid operator $\sigma\colon
\Gamma_{\inv}^{\otimes 2}\rightarrow\Gamma_{\inv}^{\otimes 2}$
associated to $\Gamma$. The algebra $\Sigma$ is always {\it commutative},
and it plays the role of the algebra of invariant polynomials over the Lie
algebra of the classical structure group. The algebras $\Sigma$ and $H(\WGi)$
are generally non-isomorphic, however they can be related \cite{d3}
through a special long exact sequence. 

\enlargethispage*{2\baselineskip}
\section{Concluding Remarks}

The spaces corresponding to algebras $\OG{d}$ are
classifying only in the relative sense, they work for bundles having a fixed
degree of complexity. The `true' quantum classifying space $\QBG{}$, together
with its universal bundle $\QEG{}$, can be constructed \cite{d4}
by taking the appropriate topological inverse limit of
algebras $\OG{d}$. However, effectively we always work with
some $\OG{d}$, for sufficiently large $d$.

Let us assume that $G$ is a classical compact Lie group.
Interestingly, the corresponding classifying space
$\QBG{}$ is still a quantum object. However, the classical classifying space
$\BG$ can be obtained by forcing the commutativity of the algebras $\OG{d}$.
Equivalently, the space $\BG$ is identified with the {\it classical part}
of $\QBG{}$. The difference in the natures of two classifying spaces
is responsible for the existence of purely quantum bundles with the
classical structure group $G$. Such bundles can be even constructed
over classical smooth manifolds (an example of a quantum line bundle
of this kind is considered in \cite{d3}). 

The presented classification theory radically differs from its classical
counterpart in the following point. In classical theory homotopic bundles over
the same manifold are isomorphic. In quantum theory this is no longer the
case, and it may happen for example that a classical bundle is homotopic
to a non-classical quantum bundle. However, it turns out \cite{d4}
that homotopic bundles give the same characteristic classes (under the
appropriate technical assumptions concerning the differential calculus). 

It is important to mention that quantum characteristic classes depend,
besides on the bundle $P$, also on specifications of corresponding
differential structures (in particular, the calculus $\Omega(M)$
on the base explicitly
depends on the calculi on $P$ and $G$). As far as classifying spaces are
concerned, the dependence is reduced to the choice of the calculus over $G$,
since it is natural to assume that the calculus on $\QEG{}$ is based on
the appropriate universal envelope. The choice of the higher-order
calculus on $G$
may have an essential influence to the algebra of universal
characteristic classes. For example, if we assume that the calculus on
$G$ is given by the universal envelope of $\cal{A}$,
the algebra $\Omega(\QBG{})$
will be acyclic. In contrast to that, if we assume that the higher-order
calculus is given by the braided \cite{w2} exterior algebra (of
the universal first-order calculus) then, generally, the cohomology of
the corresponding $\Omega(\QBG{})$ will be non-trivial.

The above described ways of defining quantum characteristic classes are not
unique. Another natural possibility is to first
factorize all differential algebras through the graded commutators
subcomplex, and then to pass to the cohomology classes. This thinking
is along the lines of cyclic homology theory \cite{c,lc}. Cyclic homology
classes of the classifying space are naturally mapped \cite{d-cycl}, via
classification maps, into cyclic homology classes of the corresponding
quantum spaces. Moreover, such a definition does not depend of differential
structures.

\enlargethispage*{2\baselineskip}


\begin{thebibliography}{10}
\bibitem{c} Connes A: {\it Noncommutative Geometry}, Academic Press (1994)
\bibitem{dr-cross} Doplicher S, Roberts J E: {\it Endomorphisms of
$C^*$-algebras, Cross Products and Duality for Compact Groups}, Ann Math,
{\bf 130}, 75--119 (1989)
\bibitem{d1} {\Dj}ur{\dj}evi\'c M:
{\it Geometry of Quantum Principal Bundles I}, Commun Math Phys
{\bf 175} (3) 457--521 (1996)
\bibitem{d2} {\Dj}ur{\dj}evi\'c M:
{\it Geometry of Quantum Principal Bundles II--Extended Version},
Preprint, Instituto de Matematicas, UNAM, M\'exico (1994)
\bibitem{d3} {\Dj}ur{\dj}evi\'c M:
{\it Characteristic Classes of Quantum Principal Bundles},
Preprint, Instituto de Matematicas, UNAM, M\'exico (1995)
\bibitem{d-bc} {\Dj}ur{\dj}evi\'c M:
{\it Quantum Principal Bundles and Their Characteristic Classes},
Lectures, Workshop on Quantum and Classical Gauge Theories, Stefan
Banach International Mathematical Center, Warsaw, Poland (May, 1995)
\bibitem{d4} {\Dj}ur{\dj}evi\'c M: {\it Geometry of Quantum Classifying
Spaces}, in preparation (1996)
\bibitem{d-cycl} {\Dj}ur{\dj}evi\'c M: {\it Cyclic Homology for Quantum
Classifying Spaces}, in preparation (1996)
\bibitem{lc} Loday J-L: {\it Cyclic Homology}, A series of Comprehensive
Studies in Mathematics {\bf 301}, Springer-Verlag (1992)
\bibitem{w1} Woronowicz S L: {\it Compact Matrix Pseudogroups}, Commun
Math Phys {\bf 111} 613--665 (1987)
\bibitem{w2} Woronowicz S L: {\it Differential Calculus on Compact
Matrix Pseudogroups $($Quantum Groups$)$}, Commun Math Phys
{\bf 122} 125--170 (1989)
\end{thebibliography}
\end{document}